\documentclass[11pt]{article}

\usepackage{amsmath}
\usepackage{amsthm}
\usepackage{amsfonts}
\usepackage{bbm}

\usepackage{xcolor}

\newcommand{\mmp}{\mathbb{P}}

\newcommand{\od}{\overset{d}{=}}

\newcommand{\me}{\mathbb{E}}

\newcommand{\mn}{\mathbb{N}}

\newcommand{\lin}{\underset{n\to\infty}{\lim}}

\DeclareMathOperator{\1}{\mathbbm{1}}

\newtheorem{thm}{Theorem}[section]

\theoremstyle{definition}

\theoremstyle{remark}
\newtheorem{rem}[thm]{Remark}

\begin{document}
\title{A remark on the paper "Renorming divergent perpetuities"}

\author{Alexander Iksanov\footnote{Faculty of Cybernetics, National Taras Shevchenko University of Kyiv, Kyiv,
Ukraine,
\newline e-mail: iksan@univ.kiev.ua} \ \ and \ \
Andrey Pilipenko\footnote{Institute of Mathematics, National
Academy of Sciences of Ukraine, Kyiv, Ukraine,
\newline e-mail: pilipenko.ay@yandex.ua} }
\maketitle
\begin{abstract}
\noindent Let $(\xi_k)$ and $(\eta_k)$ be infinite independent
samples from different distributions. We prove a functional limit
theorem for the maximum of a perturbed random walk
$\underset{0\leq k\leq n}{\max}\,(\xi_1+\ldots+\xi_k+\eta_{k+1})$
in a situation where its asymptotics is affected by both
$\underset{0\leq k\leq n}{\max}\,(\xi_1+\ldots+\xi_k)$ and
$\underset{1\leq k\leq n}{\max}\,\eta_k$ to a comparable extent.
This solves an open problem that we learned from the paper
``Renorming divergent perpetuities'' by P. Hitczenko and J.
Weso{\l}owski.
\end{abstract}

%
%
\section{Introduction and results}

Let $(\xi_k,\eta_k)_{k\in\mn}$ be a sequence of i.i.d.
two-dimensional random vectors with generic copy $(\xi,\eta)$. Let
$(S_n)_{n\in\mn_0}$ be the zero-delayed random walk with
increments $\xi_k$ for $k\in\mn$, i.e., $$S_0:=0 \ \ \text{and} \
\ S_n:=\xi_1+\ldots+\xi_n, \ \ n\in\mn.$$ Assuming that
\begin{equation}\label{10}
\me \xi=0 \ \ \text{and} \ \ v^2:={\rm Var}\,\xi<\infty,
\end{equation}
Hitczenko and Weso{\l}owski in \cite{Hitczenko+Wesolowski:2011}
investigated weak convergence of the one-dimensional distributions
of $a_n\underset{0\leq k\leq n}{\max}\,(S_k+\eta_{k+1})$ as
$n\to\infty$ for appropriate deterministic sequences $(a_n)$. More
precisely, in the proof of Theorem 3 in
\cite{{Hitczenko+Wesolowski:2011}} it is shown that (I) whenever
$\underset{0\leq k\leq n}{\max}\,S_k$ dominates $\underset{1\leq
k\leq n+1}{\max}\,\eta_k$ the limit law of $a_n\underset{0\leq
k\leq n}{\max}\,(S_k+\eta_{k+1})$ coincides with the limit law of
$a_n\underset{0\leq k\leq n}{\max}\,S_k$ which is the law of
$|B(1)|$ where $(B(t))_{t\geq 0}$ is a Brownian motion; and that
(II) whenever $\underset{1\leq k\leq n+1}{\max}\,\eta_k$ dominates
$\underset{0\leq k\leq n}{\max}\,S_k$ the limit law coincides with
that of $a_n\underset{1\leq k\leq n+1}{\max}\,\eta_k$ which is a
Fr\'{e}chet law under a regular variation assumption.

If in addition to \eqref{10} condition
\begin{equation}\label{11}
\mmp\{\eta>x\} \ \sim \ cx^{-2}, \ \ x\to\infty
\end{equation}
holds for some $c>0$, then contributions of $\underset{0\leq k\leq
n}{\max}\,S_k$ and $\underset{1\leq k\leq n+1}{\max}\,\eta_k$ to
the asymptotic behavior of $\underset{0\leq k\leq
n}{\max}\,(S_k+\eta_{k+1})$ are comparable. Hitczenko and
Weso{\l}owski conjectured (see Remark on p.~889 in
\cite{Hitczenko+Wesolowski:2011}) that whenever conditions
\eqref{10} and \eqref{11} hold, and $\xi$ and $\eta$ are
independent, the limit random variable is $\theta+vB(1)$, where
$\theta$ is independent of $B(1)$ and has a Fr\'{e}chet
distribution with parameters $2$ and $c$. Under conditions
\eqref{10} and \eqref{11} (not assuming that $\xi$ and $\eta$ are
independent) we state a functional limit result for
$n^{-1/2}\underset{0\leq k\leq [n \cdot]}{\max}\,(S_k+\eta_{k+1})$
in Theorem \ref{main} which implies that the conjecture is
erroneous (see Remark \ref{main3}).

Denote by $D:=D[0,\infty)$ the Skorokhod space of real-valued
right-continuous functions which are defined on $[0,\infty)$ and
have finite limits from the left at each positive point.
Throughout the note we assume that $D$ is equipped with the
$J_1$-topology. For $c>0$ defined in \eqref{11}, let
$N^{(c)}:=\sum_k \varepsilon_{(t_k,\,j_k)}$ be a Poisson random
measure on $[0,\infty)\times (0,\infty]$ with mean measure
$\mathbb{LEB}\times \mu_c$, where $\varepsilon_{(t,\,x)}$ is the
probability measure concentrated at $(t,x)\subset [0,\infty)\times
(0,\infty]$, $\mathbb{LEB}$ is the Lebesgue measure on
$[0,\infty)$, and $\mu_c$ is a measure on $(0,\infty]$ defined by
$$\mu_c\big((x,\infty]\big)=cx^{-2}, \ \ x>0.$$ Also, let $(B(t))_{t\geq
0}$ be a Brownian motion independent of $N^{(c)}$.
\begin{thm}\label{main}
Suppose \eqref{10} and \eqref{11}. Then
$$n^{-1/2}\underset{0\leq k\leq [n \cdot]}{\max}\,(S_k+\eta_{k+1}) \ \Rightarrow \ \ \underset{t_k\leq
\cdot}{\sup}(vB(t_k)+j_k) \ \ \text{as} \ \ n\to\infty$$ in $D$.
\end{thm}
\begin{rem}\label{main3}
Observe that $\mmp\{\underset{t_k\leq 1}{\sup}(vB(t_k)+j_k)\geq
0\}=1$, whereas $\mmp\{\theta+vB(1)<0\}>0$. This disproves the
conjecture stated in \cite{Hitczenko+Wesolowski:2011}. We note in
passing that the law of $\underset{t_k\leq 1}{\sup}(vB(t_k)+j_k)$
is different from that of $\theta+v|B(1)|\od \underset{t_k\leq
1}{\sup}j_k+ v\,\underset{t\leq 1}{\sup}B(t)$, for
$\underset{t_k\leq 1}{\sup}(vB(t_k)+j_k)<\underset{t_k\leq
1}{\sup}j_k+v\,\underset{t\leq 1}{\sup}B(t)$ a.s.
\end{rem}

After the present note was ready for submission we learned that a
version of Theorem \ref{main}, with $\xi$ and $\eta$ being
independent, has also been proved, independently and at the same
time, in \cite{Wang:2014} via a more complicated argument.

Let $C:=C[0,\infty)$ be the set of continuous functions defined on
$[0,\infty)$. Denote by $M_p$ the set of Radon point measures
$\nu$ on $[0,\infty)\times (-\infty,\infty]$ which satisfy
\begin{equation}\label{1}
\nu([0,T]\times\left\{(-\infty,-\delta]\cup[\delta,\infty]\right\})<\infty
\end{equation}
for all $\delta>0$ and all $T>0$. The $M_p$ is endowed with the
vague topology. Define the functional $F$ from $D\times M_p$ to
$D$ by
\begin{equation*}
F\left(f,\nu\right)(t):=
\begin{cases}
        \underset{k:\ \tau_k\leq t}{\sup}(f(\tau_k)+y_k),  & \text{if} \ \tau_k\leq t \ \text{for
some} \ k,\\
        f(0), & \text{otherwise},
\end{cases}
\end{equation*}
where $\nu = \sum_k \varepsilon_{(\tau_k,\,y_k)}$. Assumption
\eqref{1} ensures that $F(f,\nu)\in D$. If \eqref{1} does not
hold, $F(f,\nu)$ may lost right-continuity.
\begin{thm}\label{main2}
For $n\in\mn$, let $f_n\in D$ and $\nu_n\in M_p$. Assume that
$f_0\in C$ and that $\nu_0([0,\infty)\times(-\infty,0])=0$,
$\nu_0(\{0\}\times (-\infty,+\infty])=0$ and
$\nu_0((a,b)\times(0,\infty])\geq 1$ for all positive $a$ and $b$
such that $a<b$.

\noindent If
\begin{equation}\label{1.2}
\lin f_n= f_0 \ \ \text{in } \ D
\end{equation}
and
\begin{equation}\label{1.2'}
\lin \1_{[0,\infty)\times(0,\infty]}\nu_n=\nu_0 \mbox{ in } M_p,
\end{equation}
then
\begin{equation}\label{1.3}
\lin F(f_n,\nu_n)= F(f_0,\nu_0)
\end{equation}
in $D$.
\end{thm}
\begin{rem}
Let $a>0$ and $(T_n)_{n\in\mn_0}$ be a random sequence independent
of $(\eta_k)_{k\in\mn}$. Further, denote by $X$ a random process
with a.s. continuous paths which is independent of
$(t^\ast_k,\,j^\ast_k)$ the atoms of a Poisson random measure on
$[0,\infty)\times (0,\infty]$ with mean measure
$\mathbb{LEB}\times \mu_{c,\, a}$, where $\mu_{c,\,a}$ is a
measure on $(0,\infty]$ defined by
$\mu_{c,\,a}\big((x,\infty]\big)=cx^{-a}$, $x>0$. Whenever
\eqref{11} holds with $2$ replaced by $a$ and
$n^{-1/a}T_{[n\cdot]} \ \Rightarrow \ X(\cdot)$ in $D$, an
application of Theorem \ref{main2} allows us to infer
$$n^{-1/a}\underset{0\leq k\leq [n \cdot]}{\max}\,(T_k+\eta_{k+1}) \ \Rightarrow \ \ \underset{t^\ast_k\leq
\cdot}{\sup}(X(t^\ast_k)+j^\ast_k) \ \ \text{as} \ \ n\to\infty$$
in $D$. Details can be found in the proof of Theorem \ref{main}.
\end{rem}

The proofs of Theorem \ref{main2} and Theorem \ref{main} are given
in Section \ref{22} and Section \ref{23}, respectively.

\section{Proof of Theorem \ref{main2}}\label{22}

It suffices to prove convergence \eqref{1.3} in $D[0,T]$ for any
$T>0$ such that $\nu_0(\{T\}\times(0,\infty])=0$ (the last
condition ensures that $F(f_0,\nu_0)$ is continuous at $T$).

Let $d_T$ be the standard Skorokhod metric in $D[0,T]$. Then
$$
d_T(F(f_n,\nu_n), F(f_0,\nu_0))\leq
$$
$$
\leq d_T(F(f_n,\nu_n), F(f_0,\nu_n))+d_T(F(f_0,\nu_n), F(f_0,\nu_0))\leq
$$
$$
\leq \sup_{t\in[0,T]}|F(f_n,\nu_n)- F(f_0,\nu_n)|+d_T(F(f_0,\nu_n), F(f_0,\nu_0))\leq
$$
$$
\leq \sup_{t\in[0,T]}|f_n(t)- f_0(t)|+d_T(F(f_0,\nu_n),
F(f_0,\nu_0)),
$$
having utilized the fact that $d_T$ is dominated by the uniform
metric. It follows from \eqref{1.2} and the continuity of $f_0$
that $\lin f_n= f_0$ uniformly on $[0,T]$. Therefore we are left
with checking that
\begin{equation}\label{2.1}
\lin d_T(F(f_0,\nu_n), F(f_0,\nu_0))=0.
\end{equation}
Let $\alpha=\{0=s_0<s_1<\dots <s_m=T\}$ be a partition of $[0,T]$
such that
\begin{equation*}\label{2.2}
\nu_0(\{s_k\}\times (0,\infty])=0, \ \ k=1,...,m.
\end{equation*}
Pick now $\gamma>0$ so small that
\begin{equation*}\label{2.1'}
\nu_0((s_k,s_{k+1})\times (\gamma,\infty])\geq 1, \ \ k=0,...,m-1.
\end{equation*}

Condition \eqref{1.2'} implies that $\nu_0([0,T]\times
(\gamma,\infty])=\nu_n([0,T]\times (\gamma,\infty])=p$ for large
enough $n$ and some $p\geq 1$. Denote by
$(\bar{\tau}_i,\bar{y}_i)_{1\leq i\leq p}$ an enumeration of the
points of $\nu_0$ in $[0,T]\times (\gamma,\infty]$ with
$\bar{\tau}_1\leq \bar{\tau}_2\leq\ldots \leq \bar{\tau}_p$ and by
$(\bar{\tau}_i^{(n)}, \bar{y}_i^{(n)})_{1\leq i\leq p}$ the
analogous enumeration of the points of $\nu_n$ in $[0,T]\times
(\gamma,\infty]$. Then
\begin{equation}\label{2.2}
\lim_{n\to\infty}\sum_{i=1}^p(| \bar{\tau}^{(n)}_i- \bar{\tau}_i
|+|\bar{y}^{(n)}_i-\bar{y}_i|)=0.
\end{equation}

Define $\lambda_n$ to be continuous and strictly increasing
functions on $[0,T]$ with $\lambda_n(0) =0$, $\lambda_n(T) =T$,
$\lambda_n(\bar{\tau}^{(n)}_i)=\bar{\tau}_i$ for $i=1,\ldots,p$,
and let $\lambda_n$ be linearly interpolated elsewhere on $[0,T]$.
The relation $\lin \sup_{t\in[0,T]}|\lambda_n(t)-t|=0$ is easily
checked. Further, write
\begin{eqnarray*}
d_T(F(f_0,\nu_n), F(f_0,\nu_0))&\leq&
\sup_{t\in[0,T]}|\sup_{\lambda_n(\tau^{(n)}_k)\leq
t}(f_0(\tau^{(n)}_k)+y^{(n)}_k)-
\sup_{\lambda_n(\bar{\tau}^{(n)}_i)\leq
t}(f_0(\bar{\tau}^{(n)}_i)+\bar{y}^{(n)}_i)|\nonumber \\&+&
\sup_{t\in[0,T]}|\sup_{\tau_k\leq t}(f_0(\tau_k)+y_k)-
\sup_{\bar{\tau}_i\leq
t}(f_0(\bar{\tau}_i)+\bar{y}_i)|\\&+&\sum_{i=1}^p(|f_0(\bar{\tau}^{(n)}_i)-f_0(\bar{\tau}_i)|+|\bar{y}^{(n)}_i-\bar{y}_i|).
\end{eqnarray*}
Using \eqref{2.2} we infer
\begin{equation}\label{5}
\lin
\sum_{i=1}^p(|f_0(\bar{\tau}^{(n)}_i)-f_0(\bar{\tau}_i)|+|\bar{y}^{(n)}_i-\bar{y}_i|)=0
\end{equation}
because $f_0$ is continuous.

To proceed, put $|\alpha|:=\max_i(s_{i+1}-s_i)$ and let
$\omega_{f_0}(\varepsilon):=\underset{|u-v|<\varepsilon,\,u,v\geq
0}{\sup}\,|f_0(u)-f_0(v)|$ denote the modulus of continuity of
$f_0$. We have\footnote{We recall that
$\sup_{\lambda_n(\tau^{(n)}_k)\leq
t}(f_0(\tau^{(n)}_k)+y^{(n)}_k)=f_0(0)$ if the supremum is taken
over the empty set. Analogously for
$\sup_{\lambda_n(\bar{\tau}^{(n)}_i)\leq
t}(f_0(\bar{\tau}^{(n)}_i)+\bar{y}^{(n)}_i)$. Note that under the
assumptions of the theorem these suprema converge to $f_0(0)$ as
$t\downarrow 0$.}
\begin{eqnarray}\label{3}
&&|\sup_{\lambda_n(\tau^{(n)}_k)\leq
t}(f_0(\tau^{(n)}_k)+y^{(n)}_k)-
\sup_{\lambda_n(\bar{\tau}^{(n)}_i)\leq
t}(f_0(\bar{\tau}^{(n)}_i)+\bar{y}^{(n)}_i)|\nonumber\\&\leq&
\sup_{\lambda_n(\tau^{(n)}_k)\leq t,\
\lambda_n(\bar{\tau}^{(n)}_i)\leq t,\, \tau^{(n)}_k\neq
\bar{\tau}^{(n)}_i} \left( f_0(\tau^{(n)}_k)+y^{(n)}_k -
f_0(\bar{\tau}^{(n)}_i)\right)\vee 0 \leq
\omega_{f_0}(3|\alpha|)+\gamma.
\end{eqnarray}
Indeed, since, for $k=1,\ldots,m$,
$$(s_k, s_{k+1})\cap\{\bar{\tau}^{(n)}_1,\ldots, \bar{\tau}^{(n)}_p\}
\neq\oslash,$$ we conclude that whenever
$\lambda_n(\tau^{(n)}_k)\leq t$ there exists $\bar{\tau}^{(n)}_i$
such that $\lambda_n(\bar{\tau}^{(n)}_i)\leq t$ and
$|\tau^{(n)}_k-\bar{\tau}^{(n)}_i|\leq 3|\alpha|$. Further, all
$y_k^{(n)}$ other than $\bar{y}_i^{(n)}$, $i=1,\ldots,p$ are
smaller than or equal to $\gamma$. This explains the appearance of
the second term on the right-hand side of \eqref{3}. Arguing
similarly we infer
\begin{equation}\label{2}
|\sup_{\tau_k\leq t}(f_0(\tau_k)+y_k)- \sup_{\bar{\tau}_i\leq
t}(f_0(\bar{\tau}_i)+\bar{y}_i)|\leq
\omega_{f_0}(3|\alpha|)+\gamma.
\end{equation}
Sending in \eqref{3} and \eqref{2} $|\alpha|$ and $\gamma$ to zero
and recalling \eqref{5} we arrive at \eqref{2.1}. The proof is
complete.

\section{Proof of Theorem \ref{main}}\label{23}

According to Donsker's theorem assumption \eqref{10} implies
\begin{equation}\label{6}
n^{-1/2}S_{[n\cdot]} \Rightarrow vB(\cdot), \ \ \text{as} \ \
n\to\infty \ \ \text{in} \ \ D.
\end{equation}
It is a standard fact of the point processes theory that condition
\eqref{11} entails
\begin{equation}\label{7}
\sum_{k\geq 0}\1_{\{\eta_{k+1}>0\}}\varepsilon_{(n^{-1}k,\,
n^{-1/2}\eta_{k+1})} \ \Rightarrow \ \widehat{N}^{(c)} \ \
\text{as} \ n\to\infty \ \ \text{in} \ \ M_p,
\end{equation}
see, for instance, Corollary 4.19 (ii) in \cite{Resnick:1987}.
Here, $\widehat{N}^{(c)}$ has the same law as $N^{(c)}$ but may
depend on $B$.

In order to prove that $B$ and $\widehat{N}^{(c)}$ are actually
independent, it suffices to check that $\widehat{N}^{(c)}([0,
s]\times (\delta,\infty])$ and $B$ are independent for each fixed
$s>0$ and each fixed $\delta>0$. Fix $\delta>0$ and $s>0$ and put
$\theta^{\leq, n}_0:=0$ and $\theta^{>, n}_0:=0$, and then
$$\theta^{\leq, n}_k:=\inf\{j>\theta^{\leq, n}_{k-1}:
\eta_j\leq \sqrt{n}\delta\} \ \text{and} \ \theta^{>,
n}_k:=\inf\{j>\theta^{>, n}_{k-1}: \eta_j>\sqrt{n}\delta\}$$ for
$k\in\mn$. Further, we set $$K_n^{\leq}:=\#\{k\in\mn:
\theta_k^{\leq ,n} \leq n\} \ \ \text{and} \ \ K_n^>:=\#\{k\in\mn:
\theta_k^{> ,n} \leq n\}.$$ Then $(\xi_{\theta_k^{\leq,
n}})_{k\in\mn}$ are i.i.d. with generic copy
$\xi_{\theta^{\leq,n}}$ having the law $\mmp\{\xi_{\theta^{\leq,
n}}\in \cdot\}=\mmp\{\xi\in\cdot |\eta\leq \sqrt{n}\delta\}$,
while $(\xi_{\theta_k^{>, n}})_{k\in\mn}$ are i.i.d. with generic
copy $\xi_{\theta^{>,n}}$ having the law $\mmp\{\xi_{\theta^{>,
n}}\in \cdot\}=\mmp\{\xi\in \cdot|\eta>\sqrt{n}\delta\}$. For any
$\varepsilon>0$,
$$\mmp\{|\xi_{\theta^{>, n}}|>\sqrt{n}\varepsilon\}\leq
\mmp\{|\xi_{\theta^{>,
n}}|>\sqrt{n}\varepsilon\}/\mmp\{\eta>\sqrt{n}\delta\}\sim
c^{-1}\delta^2 n\mmp\{|\xi_{\theta^{>,
n}}|>\sqrt{n}\varepsilon\}$$ which proves that $\lin
n^{-1/2}\xi_{\theta^{>, n}}=0$ in probability. Since
$K^>_{[nT]}=\sum_{k\geq 0}\varepsilon_{(n^{-1}k,\,
n^{-1/2}\eta_{k+1})}([0,T]\times (\delta,\infty])$ converges to
$\widehat{N}^{(c)}([0,T]\times (\delta,\infty])$ in distribution
as $n\to\infty$, the right-hand side of
\begin{equation*}\label{4}
n^{-1/2}\underset{t\in
[0,\,T]}{\sup}\bigg|\sum_{i=1}^{[nt]}\xi_i-\sum_{j=1}^{K^\leq_{[nt]}}\xi_{\theta^{\leq,n}_j}\bigg|=
n^{-1/2}\underset{t\in
[0,\,T]}{\sup}\bigg|\sum_{k=1}^{K^>_{[nt]}}\xi_{\theta^{>,n}_k}\bigg|\leq
n^{-1/2}\sum_{k=1}^{K^>_{[nT]}}|\xi_{\theta^{>,n}_k}|,
\end{equation*}
where $T>0$ is arbitrary, converges to zero in probability, as
$n\to\infty$. Therefore,
$$n^{-1/2}\sum_{j=1}^{K^\leq_{[n\cdot]}}\xi_{\theta^{\leq,n}_j} \
\Rightarrow \ vB(\cdot) \ \ \text{as} \ \ n\to\infty \ \ \text{in}
\ \ D.$$ Observe further that $n^{-1}K^\leq_{[n\cdot]}\Rightarrow
f(\cdot)$ as $n\to\infty$ in $D$, where $f(t)=t$, $t\geq 0$ which
implies
$$\bigg(n^{-1/2}\sum_{j=1}^{K^\leq_{[n\cdot]}}\xi_{\theta^{\leq,n}_j}, n^{-1/2}\sum_{j=1}^{[n\cdot]}\xi_{\theta^{\leq,n}_j}\bigg)
\ \Rightarrow \ \big(vB(\cdot), vB(\cdot)\big) \ \ \text{as} \
n\to\infty \ \ \text{in} \ D\times D.$$ Since $K^>_{[ns]}$ is
independent of $(\xi_{\theta^{\leq, n}_k})_{k\in\mn}$ we conclude
that $B$ and $\widehat{N}^{(c)}([0,s]\times (\delta,\infty])$ are
independent, as claimed.

Using the independence of $B$ and $\widehat{N}^{(c)}$, relations
\eqref{6} and \eqref{7} can be combined into the joint convergence
\begin{equation*}\label{8}
\bigg(n^{-1/2}S_{[n\cdot]},\sum_{k\geq
0}\1_{\{\eta_{k+1}>0\}}\varepsilon_{(n^{-1}k,\,
n^{-1/2}\eta_{k+1})}\bigg) \ \Rightarrow \ \big(vB(\cdot),
\widehat{N}^{(c)}\big) \ \text{as} \ \ n\to\infty \ \text{in} \ \
D\times M_p
\end{equation*}
(in the product topology). By the Skorokhod representation theorem
there are versions which converge a.s. Retaining the original
notation for these versions we apply Theorem \ref{main2} with
$f_n(\cdot)=n^{-1/2}S_{[n\cdot]}$, $f_0=vB$, $\nu_n=\sum_{k\geq
0}\varepsilon_{\{n^{-1}k,\, n^{-1/2}\eta_{k+1}\}}$ and $\nu_0=
\widehat{N}^{(c)}$. We already know that conditions \eqref{1.2}
and \eqref{1.2'} are fulfilled. Plainly,
$\widehat{N}^{(c)}([0,\infty)\times(-\infty,0])=0$ a.s. and
$\widehat{N}^{(c)}(\{0\}\times (-\infty,+\infty])=0$ a.s. Further
$\widehat{N}^{(c)}([0,T]\times\left\{(-\infty,-\delta]\cup[\delta,\infty]\right\})<\infty$
a.s. for all $\delta>0$ and all $T>0$ because
$\mu_c((-\infty,-\delta]\cup[\delta,\infty])=
\mu_c([\delta,\infty]) <\infty$, and
$\widehat{N}^{(c)}((a,b)\times(0,\infty])\geq 1$ a.s. whenever
$0<a<b$ because  $\mu_c((0,\infty])=\infty$. The proof is
complete.

\end{document}